\documentclass[12pt]{article}
\usepackage{enumerate}
\usepackage{amsfonts}
\usepackage{amsmath}
\usepackage{amssymb}

\newcounter{defin}  \newcounter{lemma}  \newcounter{theorem}
\newcounter{property} \newcounter{corol}  \newcounter{remark} \newcounter{example}

\begin{document}
\title{Characterization of convex $\mu$-compact sets}

\author{M.E.~Shirokov\thanks{e-mail:msh@mi.ras.ru}\\\\
Steklov Mathematical Institute, Moscow, Russia}

\date{} \maketitle

\begin{abstract}
The class of $\mu$-compact sets can be considered as a natural
extension of the class of compact metrizable subsets of locally
convex spaces, to which the particular results well known for
compact sets can be generalized. This class contains all compact
sets as well as many noncompact sets widely used in applications. In
this paper we give a characterization of a convex $\mu$-compact set
in terms of properties of functions defined on this set. Namely, we
prove that the class of convex $\mu$-compact sets can be
characterized by continuity of the operation of convex closure of a
function (= the double Fenchel transform) with respect to monotonic
pointwise converging sequences of continuous bounded and of lower
semicontinuous lower bounded functions.
\end{abstract}\bigskip\bigskip\bigskip

The properties of compact sets in the context of convex analysis
have been studied by many authors (see~\cite{1,4,5} and the
references therein). It is natural to ask about possible
generalizations of the results proved for compact convex sets to
noncompact sets. In ~\cite{2} one such generalization concerning the
particular class of sets called $\mu$-compact sets is considered.
In~\cite{2,3} it is shown that for this class of sets, which
includes all compact convex sets as well as some noncompact sets
widely used in applications, many results of the Choquet
theory~\cite{1} and~of the Vesterstrom-O'Brien theory ~\cite{4,5}
can be proved. In this paper we give a characterization of a convex
$\mu$-compact set in terms of properties of functions defined on
this set.\medskip

In what follows $\mathcal{A}$~is a bounded convex complete separable
metrizable subset of some locally convex space.\footnote{This means
that the topology on the set $\mathcal{A}$ is defined by a countable
subset of the family of seminorms, generating the topology of the
entire locally convex space, and this set is separable and complete
in the metric generated by this subset of seminorms.} Let
$C(\mathcal{A})$~be the set of all continuous bounded functions on
the set $\mathcal{A}$ and ~$M(\mathcal{A})$~be the set of all Borel
probability measures on the set $\mathcal{A}$ endowed with the weak
convergence topology~\cite[Chapter II, \S 6]{6}. Let $\mathrm{co}f$
and $\overline{\mathrm{co}}f$ be the convex hull and the convex
closure of a function~$f$, which are defined respectively as the
maximal convex and the maximal convex closed (that is, lower
semicontinuous) functions majorized by $f$ \cite{7}.\footnote{The
convex closure of a function is also called the lower (convex)
envelope of this function \cite{1}.}\smallskip

With an arbitrary measure $\mu\in M(\mathcal{A})$ we associate its
barycenter (average)~$\mathbf b\, (\mu)\, \in \, \mathcal{A}$, which
is defined by the Pettis integral~(see \cite{V&T,12})
\begin{equation}\label{eq1}
\mathbf b(\mu)\ =\ \int_{\mathcal{A}}x \mu(dx).
\end{equation}
For arbitrary $x\in\mathcal{A}\,$ let $M_{x}(\mathcal{A})$ be a
convex closed subset of the set $M(\mathcal{A})$ consisting of such
measures  $\mu$  that $\,\mathbf b(\mu) = x$.\smallskip

The barycenter map
\begin{equation}\label{b-map}
M(\mathcal{A})\ \ni \ \mu\; \mapsto \; \mathbf b(\mu)\ \in
 \mathcal{A}
\end{equation}
is continuous (this can be shown easily by applying Prokhorov's
theorem~\cite[Ch.II, Th.6.7]{6}). Hence the image of any compact
subset of $M(\mathcal{A})$ under this  map is a compact subset of
$\mathcal{A}$. The $\mu$-compact sets are defined in \cite{2,3} by
the converse requirement.\medskip

\textbf{Definition.} A set $\mathcal{A}$ is called $\,\mu$-{\it
compact } if the preimage of any compact subset of $\mathcal{A}$
under barycenter map (\ref{b-map}) is a compact subset of
$M(\mathcal{A})$.
\medskip

Any compact set is $\mu$-compact, since compactness of $\mathcal{A}$
implies compactness of $M(\mathcal{A})$~\cite{6}. The
$\mu$-compactness property is studied in detail in~\cite{3}, where
simple criteria of this property have been established. By using
these criteria  $\mu$-compactness of the following noncompact sets
has been proved:
\begin{itemize}

\item[--]
the positive parts of the unit balls of the Banach space $\ell_{1}$
and of the Banach space $\mathfrak{T}(\mathcal{H})$ of trace class
operators in a separable Hilbert space $\mathcal{H}$;

\item[--]
the set of positive Borel measures on an \textit{arbitrary} complete
separable metric space with the total variation $\leq 1$ endowed
with the weak convergence topology;

\item[--]
the positive parts of the unit balls of the Banach spaces of linear
bounded operators in $\ell_{1}$ and in~$\mathfrak{T}(\mathcal{H})$
endowed with the strong operator topology.
\end{itemize}
In particular, this implies $\mu$-compactness of the set of all
Borel probability measures on an arbitrary complete separable metric
space endowed with the weak convergence topology, of the set of
quantum states and of the set of quantum operations endowed with the
strong operator topology~\cite{8}.\smallskip

It is essential to note that the $\mu$-compactness property of a
convex set is not purely topological but reflects the special
relation between the topology and the convex structure of this
set~\cite{3}.\smallskip

The following theorem shows that the class of convex $\mu$-compact
sets can be characterized by continuity of the operation of convex
closure (coinciding with the double Fenchel transform) with respect
to monotonic pointwise converging sequences of functions.\pagebreak

\textbf{Theorem.} \emph{The following properties are equivalent:}
\begin{enumerate}[\rm(i)]

\item
\emph{the set $\mathcal{A}$ is $\mu$-compact;}

\item
\emph{for an arbitrary increasing sequence $\{f_{n}\}\subset
C(\mathcal{A})$, converging pointwise to a function $f_{0}\in
C(\mathcal{A})$, the sequence $\{\overline{\mathrm{co}}f_{n}\}$
converges pointwise to the function
$\,\overline{\mathrm{co}}f_{0}$;}

\item
\emph{for an arbitrary increasing sequence  $\{f_{n}\}$ of lower
semicontinuous lower bounded functions on $\mathcal{A}$, converging
pointwise to a function $f_{0}$, the sequence
$\{\overline{\mathrm{co}}f_{n}\}$  converges pointwise to the
function $\,\overline{\mathrm{co}}f_{0}$.}
\end{enumerate}
\emph{If these equivalent properties hold then for an arbitrary
decreasing sequence $\{f_{n}\}$ of lower semicontinuous bounded
functions on~the set $\mathcal{A}$,  converging pointwise to a lower
semicontinuous bounded function $f_{0}$, the sequence
$\{\overline{\mathrm{co}}f_{n}\}$  converges pointwise to the
function $\overline{\mathrm{co}}f_{0}$.} \smallskip

\textbf{Remark 1.} The functions $\,\overline{\mathrm{co}}f_{n}\,$
in (ii) are not necessarily continuous (only lower semicontinuous).
By the generalized Vesterstrom-O'Brien theorem (Theorem 1 in
\cite{3}) these functions are continuous provided the set
$\mathcal{A}$ is $\mu$-compact and \emph{stable} (the last property
means openness of the convex mixing map
$\mathcal{A}\times\mathcal{A}\ni(x,y)\mapsto \frac{1}{2}(x + y) \in
\mathcal{A}$ \cite{11}). \medskip

\textbf{Proof.} It is sufficient to show that
$\mathrm{(i)}\Rightarrow\mathrm{(iii)}$ and
$\mathrm{(ii)}\Rightarrow\mathrm{(i)}$.\smallskip

$\mathrm{(i)}\Rightarrow\mathrm{(iii)}\,$ Let $\{f_{n}\}$ be an
increasing sequence of lower semicontinuous lower bounded  functions
on the set $\mathcal{A}$, converging pointwise to a function
$f_{0}$. We may assume that this sequence consists of nonnegative
functions. It suffices to show that the assumption on existence of
such $x_{0}\in\mathcal{A}$ that
$$
\overline{\mathrm{co}}f_{n}(x_{0})\leq\overline{\mathrm{co}}f_{0}(x_{0})-\Delta,
\quad \Delta>0,\quad \forall n,
$$
where "$\leq+\infty-\Delta$" means "$\leq\Delta$", leads to a
contradiction.\smallskip

By Proposition 6 in \cite{3} we have
\begin{equation}\label{co-rep}
\overline{\mathrm{co}}f_{n}(x_{0})=\inf_{\mu\in
M_{x_{0}}(\mathcal{A})}\mu(f_{n}),\quad n=0,1,2...,\quad
\textup{where}\;\, \mu(f)=\int_{\mathcal{A}} f(y)\,\mu(dy),
\end{equation}
and this infimum is attained at a particular measure $\mu_{n}$ in
$M_{x_{0}}(\mathcal{A})$, t.i.
$\overline{\mathrm{co}}f_{n}(x_{0})=\mu_{n}(f_{n})$.

For definiteness suppose
$\,\overline{\mathrm{co}}f_{0}(x_{0})\!<\!+\infty$ (the case
$\,\overline{\mathrm{co}}f_{0}(x_{0})\!=\!+\infty\,$ is considered
similarly). By Fenchel's theorem (see \cite{7}) there exists a
continuous affine function $\,\alpha\,$ on $\mathcal{A}$ such that
\begin{equation}\label{ineq-p-1}
\alpha(x)\leq f_{0}(x),\quad\forall
x\in\mathcal{A},\quad\textup{ш}\quad\overline{\mathrm{co}}f_{0}(x_{0})\leq\alpha(x_{0})+\textstyle{\frac{1}{2}}\Delta.
\end{equation}
Since the function $\alpha$ is affine, we have
\begin{equation}\label{ineq-p-2}
\begin{array}{c}
\displaystyle
\mu_{n}(\alpha)-\mu_{n}(f_{n})=\alpha(x_{0})-\overline{\mathrm{co}}f_{n}(x_{0})\\\\\displaystyle=
[\alpha(x_{0})-\overline{\mathrm{co}}f_{0}(x_{0})]+[\overline{\mathrm{co}}f_{0}(x_{0})
-\overline{\mathrm{co}}f_{n}(x_{0})]\geq
-\textstyle{\frac{1}{2}}\Delta+\Delta=\frac{1}{2}\Delta.
\end{array}
\end{equation}
The assumed $\mu$-compactness of the set $\mathcal{A}$ implies
relative compactness of the sequence $\{\mu_{n}\}$. By Prokhorov's
theorem this sequence is \textit{tight}, which means that for any
$\varepsilon>0$ there exists such compact set
$\mathcal{K}_{\varepsilon}\subset\mathcal{A}$ that
$\mu(\mathcal{A}\setminus\mathcal{K}_{\varepsilon})<\varepsilon\,$
\cite{6}. Let $M=\sup_{x\in\mathcal{A}}|\alpha(x)|$ and
$\varepsilon_{0}=\frac{1}{4M}\Delta$. By using (\ref{ineq-p-2}) we
obtain
$$
\int_{\mathcal{K}_{\varepsilon_{0}}}(\alpha(x)-f_{n}(x))\mu_{n}(dx)\geq
\textstyle{\frac{1}{2}}\Delta-\displaystyle\int_{\mathcal{A}\setminus\mathcal{K}_{\varepsilon_{0}}}(\alpha(x)-f_{n}(x))\mu_{n}(dx)
\geq \textstyle{\frac{1}{4}}\Delta.
$$
Hence the set
$\mathcal{C}_{n}=\{x\in\mathcal{K}_{\varepsilon_{0}}\,|\,
\alpha(x)\geq f_{n}(x)+\frac{1}{4}\Delta\}$ is not empty for all
$n$.\smallskip

Since the sequence $\{f_{n}\}$ is increasing, the sequence
$\{\mathcal{C}_{n}\}$ of \textit{closed} subsets of the
\textit{compact} set $\mathcal{K}_{\varepsilon_{0}}$ is monotone:
$\mathcal{C}_{n+1}\subseteq\mathcal{C}_{n},\;\forall n$. Hence these
exists $x_{*}\in\mathcal{C}_{n}$ for all $n$. This means that
$\alpha(x_{*})\geq f_{n}(x_{*})+\frac{1}{4}\Delta\,$ for all  $n$
and hence $\alpha(x_{*})> f_{0}(x_{*})$ contradicting
(\ref{ineq-p-1}).
\medskip

$\mathrm{(ii)}\Rightarrow\mathrm{(i)}\,$ Suppose the set
$\mathcal{A}$ is not $\mu$-compact. Then then there exists a
sequence $\{\mu_{k}\}\in M(\mathcal{A})$, which is a not relatively
compact and such that the sequence $\{x_{k}=\mathbf{b}(\mu_{k})\}$
converges. By the below Lemma 1 one can consider that this sequence
consists of finitely supported measures. By Prokhorov's theorem the
sequence $\{\mu_{k}\}$ is not tight. The below Lemma 2 (with the
remark after it) guarantees existence of such $\varepsilon>0$ and
$\delta>0$ that for any compact set
$\mathcal{K}\subseteq\mathcal{A}$ and any natural $N$ there is such
$k>N$ that $\mu_{k}(U_{\delta}(\mathcal{K}))< 1-\varepsilon$, where
$U_{\delta}(\mathcal{K})$ is the closed $\delta$-vicinity of the set
$\mathcal{K}$ (as a subset of the metric space $\mathcal{A}$). Let
$\{\mathcal{K}_{n}\}$ be an increasing sequence of compact subsets
of $\mathcal{A}$ such that
$\bigcup_{n\in\mathbb{N}}U_{\delta/2}(\mathcal{K}_{n})=\mathcal{A}$.
Denote by $d(\cdot\,,\cdot)$ the metric in $\mathcal{A}$. For each
natural $n$ consider the continuous bounded function
\begin{equation}\label{f-n}
f_{n}(x)=1-2\,\delta^{-1}\inf_{y\in
U_{\delta/2}(\mathcal{K}_{n})}d(x,y)
\end{equation}
on the set $\mathcal{A}$ such that $f_{n}(x)=1$, if $x\in
U_{\delta/2}(\mathcal{K}_{n}),$ and $f_{n}(x)<0$, if $x\in
\mathcal{A}\setminus U_{\delta}(\mathcal{K}_{n})$. It is clear that
$f_{0}(x)\doteq\lim_{n\rightarrow+\infty}f_{n}(x)\equiv 1$ and hence
$\overline{\mathrm{co}}f_{0}(x)\equiv 1$. Let $x_{0}$ be a limit of
the sequence $\{x_{k}\}$. To obtain a contradiction it suffices to
show that
\begin{equation}\label{ineq}
\overline{\mathrm{co}}f_{n}(x_{0})<1-\varepsilon\quad \forall
n\in\mathbb{N}.
\end{equation}
By the above-stated property of the sequence $\{\mu_{k}\}$ for each
$n$ and any natural $N$ there exists such $k>N$ that
$\mu_{k}(U_{\delta}(\mathcal{K}_{n}))< 1-\varepsilon$ and hence
$$
\overline{\mathrm{co}}f_{n}(x_{k})\leq
\mathrm{co}f_{n}(x_{k})\leq\int_{\mathcal{A}}
f_{n}(x)\mu_{k}(dx)<1-\varepsilon,
$$
since $\mu_{k}$ is a measure with finite support. This inequality
and lower semicontinuity of the function
$\overline{\mathrm{co}}f_{n}$ imply (\ref{ineq}).\medskip

By using representation (\ref{co-rep}) and the monotonic convergence
theorem, it is easy to prove the last assertion of the theorem.
$\square$\medskip

In the above proof the following assertions (well known in the
measure theory) are used.\medskip

\textbf{Lemma 1.} \textit{For an arbitrary sequence
$\{\mu_{k}\}\subset M(\mathcal{A})$, which is not relatively
compact, there exists a sequence $\{\tilde{\mu}_{k}\}\subset
M(\mathcal{A})$ of finitely supported measures, which is not
relatively compact as well, such that
$\,\mathbf{b}(\tilde{\mu}_{k})=\mathbf{b}(\mu_{k})$ for all
$\,k$.}\medskip\\
Since the set $M(\mathcal{A})$ can be considered as a complete
separable metric space \cite[Chapter II]{6}, the above lemma is
easily proved by using density of the set of finitely supported
measures in $M_{x}(\mathcal{A})$ for each $x\in\mathcal{A}$
\cite[Lemma 1]{2}.
\medskip

\textbf{Lemma 2.} \textit{A subset $M_{0}\subseteq M(\mathcal{A})$
is tight if and only if for any $\varepsilon>0$ and  $\delta>0$
there exists a compact subset
$\,\mathcal{K}(\varepsilon,\delta)\subseteq\mathcal{A}$ such that
$$
\mu(U_{\delta}(\mathcal{K}(\varepsilon,\delta)))\geq 1-\varepsilon
$$
for all $\mu\in M_{0}$, where
$U_{\delta}(\mathcal{K}(\varepsilon,\delta))$ is the closed
$\delta$-vicinity of the set
$\,\mathcal{K}(\varepsilon,\delta)$.}\medskip

Since any finite subset of $M(\mathcal{A})$ is tight, the words
"\emph{for all $\mu\in M_{0}$}" in the above criterion may be
replaced by "\emph{for all $\mu\in M_{0}\setminus M'$, where $M'$ is
a finite subset of $M_{0}$}".\medskip

\textbf{Proof.} It is easy to see that tightness of the set $M_{0}$
implies validity of the condition in the lemma. Suppose this
condition holds. For arbitrary $\varepsilon>0$ and each
$n\in\mathbb{N}$ let $\mathcal{K}_{n}=\mathcal{K}\left(\varepsilon
2^{-n}, \varepsilon 2^{-n}\right)$. Then for the compact set
$\mathcal{K}=\bigcap_{n\in\mathbb{N}}U_{\varepsilon
2^{-n}}(\mathcal{K}_{n})$ we have
$$
\mu(\mathcal{A}\setminus\mathcal{K})\leq\sum_{n=1}^{+\infty}\mu(\mathcal{A}\setminus
U_{\varepsilon 2^{-n}}(\mathcal{K}_{n}))
<\sum_{n=1}^{+\infty}\varepsilon 2^{-n}<\varepsilon
$$
for all $\mu\in M_{0}$, which means that the set $M_{0}$ is tight.
$\square$\medskip

It is well known that an arbitrary lower semicontinuous lower
bounded function on a metric space can be represented as a pointwise
limit of some increasing sequence of continuous bounded functions
\cite{12}. By using the above Theorem and the generalized
Vesterstrom-O'Brien theorem mentioned in Remark 1 this observation
can be strengthened as follows. \vspace{5pt}

\textbf{Corollary.} \emph{An arbitrary lower semicontinuous lower
bounded convex (correspondingly, concave) function on a stable
$\mu$-compact set $\mathcal{A}$ can be represented as a pointwise
limit of some increasing sequence of convex (correspondingly,
concave) continuous bounded functions.}\bigskip

\textbf{Remark 2.} In \cite{3} the weaker version of the
$\mu$-compactness property of a set $\mathcal{A}$ defined by the
requirement of compactness of the set $M_{x}(\mathcal{A})$ for each
$x$ in $\mathcal{A}$ is considered. This property called
\emph{pointwise $\mu$-compactness} was used to show that even slight
relaxing of the $\mu$\nobreakdash-\hspace{0pt}compactness condition
in the generalized Vesterstrom-O'Brien theorem leads to breaking its
validity. The class of pointwise
$\mu$\nobreakdash-\hspace{0pt}compact sets is wider than the class
of $\mu$\nobreakdash-\hspace{0pt}compact sets, in particular, it
contains the simplex
$\left\{\{x_{i}\}_{i=1}^{+\infty}\,|\,x_{i}\geq0,\,\forall i,\,
\sum_{i=1}^{+\infty}x_{i}\leq 1\right\}\subset\ell_{p}$ for any
$p\geq1$, which is $\mu$\nobreakdash-\hspace{0pt}compact only for
$p=1$ \cite[Proposition 13]{3}. For an arbitrary convex pointwise
$\mu$\nobreakdash-\hspace{0pt}compact set the assertions of the
Krein-Milman theorem and of the Choquet theorem are valid
\cite[Proposition 5]{3}, but representation (\ref{co-rep}) for the
convex closure of a lower semicontinuous lower bounded function $f$
does not hold in general \cite[Example 1]{3}.\pagebreak

Similar to a convex $\mu$-compact set a convex pointwise
$\mu$-compact set can be characterized in terms of properties of
functions defined on this set. Namely, one can show that \emph{the
following properties are equivalent:}
\begin{enumerate}[\rm(i)]
\item
\emph{the set $\mathcal{A}$ is pointwise $\mu$-compact;}
\item
\emph{for an arbitrary increasing sequence $\{f_{n}\}\subset
C(\mathcal{A})$, converging pointwise to a function $f_{0}\in
C(\mathcal{A})$, the sequence $\{\mathrm{co}f_{n}\}$ converges
pointwise to the function $\,\mathrm{co}f_{0}$.}
\end{enumerate}\medskip
Since
$$
\mathrm{co}f_{n}(x)=\inf_{\mu\in
M^{f}_{x}(\mathcal{A})}\int_{\mathcal{A}} f_{n}(y) \,\mu(dy),\quad
x\in\mathcal{A},\quad n=0,1,2...,
$$
where $M^{f}_{x}(\mathcal{A})$ is a subset of $M_{x}(\mathcal{A})$
consisting of finitely supported measures, the implication
$\,\mathrm{(i)\Rightarrow(ii)}\,$ in the above assertion is proved
by noting that pointwise $\mu$\nobreakdash-\hspace{0pt}compactness
of $\mathcal{A}$ and Prokhorov's theorem implies tightness of
$\,M^{f}_{x}(\mathcal{A})\,$ and by using Dini's lemma. The
implication $\,\mathrm{(ii)\Rightarrow(i)}\,$ can be established by
using the proof of the implication $\,\mathrm{(ii)\Rightarrow(i)}\,$
in the Theorem with $x_{k}=x_{0}$ for all $k$.\medskip

The above characterization of a convex pointwise $\mu$-compact set,
the Theorem and Corollary 2 in \cite{3} show that
$$
\{\textrm{pointwise}\; \mu\textup{-compactness\,of}\;\mathcal{A}
\}\wedge\{\mathrm{co}f=\overline{\mathrm{co}}f\;\forall f\in
C(\mathcal{A})\}=\{\mu\textup{-compactness\,of}\;\mathcal{A}\}.
$$

\bigskip

The implication $\,\mathrm{(i)\Rightarrow(iii)}\,$ in the Theorem
can be used in study of entropic characteristics of quantum states
\cite[\S 6.2]{10}.

\end{document}